
\input amstex
\NoBlackBoxes
\documentstyle{amsppt}
\magnification 1100
\pageheight {8.8 truein} 
\pagewidth {6.5truein} 
\baselineskip=16pt
\TagsOnLeft

\loadeufm

\def\tarrow{\hbox{\lower1.0ex\hbox{$ 
  \buildrel{\lower0.6ex\hbox{$\textstyle\a$}}\over
{\buildrel{\lower0.5ex\hbox{$\textstyle\a$}}\over{\textstyle\a}}$}
}}

\def\darrow{\hbox{\lower0.5ex\hbox{$ 
  \buildrel{\lower0.7ex\hbox{$\textstyle\a$}}\over{\textstyle\a}$}}}

\def\dashtarrow{\hbox{\lower1.0ex\hbox{$ 
  \buildrel{\lower0.6ex\hbox{$\textstyle\a$}}\over
{\buildrel{\lower0.5ex\hbox{$\textstyle\a$}}\over{\textstyle\dashrightarrow}}$}
}}

\def\dashdarrow{\hbox{\lower0.5ex\hbox{$ 
  \buildrel{\lower0.7ex\hbox{$\textstyle\a$}}\over{\textstyle\dashrightarrow}$}}}

\def\lrdarrow{\hbox{\lower0.5ex\hbox{$ 
  \buildrel{\lower0.7ex\hbox{$\textstyle\rightarrow$}}\over{\textstyle\leftarrow}$}}}

\def \bsno{\bigskip\noindent}
\def \msno{\medskip\noindent}

\def \a{\longrightarrow}

\def \fm{\frak m}

\def \Rees{\operatorname{Rees}\;}
\def \G{\operatorname{Gr}\;}
\def \charr{\operatorname{char}\;}
\def \GrAut{\operatorname{GrAut}}

\def \Aut{\operatorname{Aut}}

\def \dim{{\operatorname{dim}}\,}
\def \Kdim{{\operatorname{Kdim}}\,}
\def \Cdim{{\operatorname{Cdim}}\,}

\def \Ext{\operatorname{Ext}}

\def \Tor{\operatorname{Tor}}

\def \Int{\operatorname{Int}\; }

\def \ad{\operatorname{ad}\; }

\def \uExt{\operatorname{\underline{Ext}}}

\def \id{\operatorname{injdim}\;}

\def \hdet{\operatorname{hdet}}
\def \det{\operatorname{det}}

\def\lrarrow{\hbox{\lower0.5ex\hbox{$ 
  \buildrel{\lower0.7ex\hbox{$\textstyle\longrightarrow$}}
     \over{\textstyle\longleftarrow}$}}}
  
\def\gldim{\operatorname{gldim}}

\def \fm{\frak m}

\def \a{\longrightarrow}

\def \GKdim{\operatorname{GKdim}}

\def\FilAut{\operatorname{FilAut}}

\document  

\topmatter
\title Gorensteinness of invariant subrings of quantum
algebras
\endtitle

\author Naihuan Jing and James J. Zhang
\endauthor

\address
(Jing) Department of Mathematics, North Carolina State University, 
Raleigh, NC 27695-8205. 
\endaddress

\email jing\@math.ncsu.edu
\endemail

\address
(Zhang) Department of Mathematics, Box 354350, University of 
Washington, Seattle,
WA 98195.  
\endaddress

\email zhang\@math.washington.edu
\endemail

\abstract
We prove Auslander-Gorenstein and $\GKdim$-Macaulay properties for
certain invariant subrings of some quantum algebras, the Weyl algebras,
and the universal enveloping algebras of finite dimensional Lie algebras.
\endabstract

\keywords group actions, invariant subrings, Auslander-Gorenstein property, 
Macaulay property
property
\endkeywords
\subjclass 16E10, 16W20
\endsubjclass
\endtopmatter
\centerline{\bf 0. Introduction} 

\bsno

Given a noncommutative algebra it is generally 
difficult to determine its homological properties
such as global dimension and injective dimension. In this
paper we use the noncommutative version of Watanabe 
theorem proved in [{\bf JoZ}, 3.3] to give some simple 
sufficient conditions for certain classes of invariant
rings having some good homological properties.

Let $k$ be a base field. Vector spaces, algebras, etc. are over $k$.
Suppose $G$ is a finite group of automorphisms of an algebra $A$. Then
the {\it invariant subring} is defined to be
$$A^G=\{x\in A\; |\; g  (x)=x \quad{\text{for all}}\; g  \in G\}.$$

Let $A$ be a filtered ring with a filtration 
$\{F_i\;|\; i \geq 0\}$ such that $F_0=k$. 
The associated graded ring is defined to be 
$\G A=\bigoplus_n F_n/F_{n-1}$. A {\it filtered (or graded)
automorphism} of $A$ (or $\G A$) is an automorphism preserving 
the filtration (or the grading). The following is a noncommutative 
version of [{\bf Ben}, 4.6.2].

\proclaim{Theorem 0.1} Suppose $A$ is a filtered ring such that the associated
graded ring $\G A$ is isomorphic to a skew polynomial 
ring $k_{p_{ij}}[x_1,\cdots, x_n]$,
where $p_{ij}\neq p_{kl}$ for all $(i,j)\neq (k,l)$.
Let $G$ be a finite
group of filtered automorphisms of $A$ with $|G|\neq 0$ in $k$. If $\det \;
g|_{(\oplus_{i=1}^n kx_i)}=1$ for all $g\in G$, then 
$A^G$ is Auslander-Gorenstein and $\GKdim$-Macaulay.
\endproclaim

The condition on $p_{ij}$ in Theorem 0.1 is needed as shown
by Example 2.9.

The definitions of Auslander-Gorenstein and 
$\GKdim$-Macaulay are given in Definition 0.6.
The skew polynomial algebra 
$k_{p_{ij}}[x_1,\cdots, x_n]$ is generated by 
$\{x_1,\cdots, x_n\}$ with relations $x_ix_j=p_{ij}x_jx_i$
for all $i<j$, where the parameters $\{p_{ij}\;|\; i<j\}$
are nonzero elements in $k$. We let $p_{ii}=1$ and $p_{ij}=p_{ji}^{-1}$, then
the relations are $x_ix_j=p_{ij}x_jx_i$.
We always think
$ k_{p_{ij}}[x_1,\cdots, x_n]$ as a connected graded ring
with the degree of $x_i$ not necessarilly $1$. We identify the graded
vector space $V=\bigoplus_{i=1}^n kx_i$ with the quotient 
space ${\frak m}/{\frak m}^2$, where $\frak m$ is 
the maximal graded ideal of  $k_{p_{ij}}[x_1,\cdots,x_n]$.
Any filtered or graded automorphism $g$ induces a linear
automorphism of $V$ and $\det g|_V$ is the usual determinant of
the linear map $g:V\to V$. The 
hypothesis on $\det g$ is necessary as shown in [Example 3.6].
 
A general version of Theorem 0.1 
holds when we replace $\det g$ by the homological
determinant [Theorems 3.3 and 3.5]. The notion of homological determinant 
was introduced in [{\bf JoZ}, 2.3] (see Section 2). 

We can extend the action of $g$ from $V$ 
to the  exterior algebra 
$\Lambda(V)$. Therefore $g$ acts on the tensor product 
$A\otimes \Lambda(V)$ if $A$ is as in Theorem 0.1. 
The following is a noncommutative version of [{\bf Ben}, 5.3.2].

\proclaim{Theorem 0.2} Suppose $A$ is a filtered ring such that the associated
graded ring $\G A$ is isomorphic to the skew polynomial ring 
$k_{p_{ij}}[x_1,\cdots, x_n]$, where $p_{ij}\neq p_{kl}$ for
$(i,j)\neq (k, l)$. Let $G$ be a finite
group of filtered automorphisms of $A$ with $|G|\neq 0$ in $k$. Then 
$(A\otimes \Lambda(V))^G$ is Auslander-Gorenstein and $\GKdim$-Macaulay 
where $V=\bigoplus_{i=1}^n kx_i$.
\endproclaim 

Theorem 3.8 is another version of Theorem 0.2, which 
contains the commutative case.

Note that in Theorem 0.2 there is no hypothesis on 
$\det \; g$ as in the commutative case [{\bf Ben}, 5.3.2]. 
In general, $(A\otimes \Lambda(V))^G$ is not prime and
hence is not Auslander-regular. 

Theorems 0.1 and 0.2 also hold if $\G A$ is  
the 4-dimensional Sklyanin algebra (by Example 2.2 and
Theorem 3.5) or the multiparameter 
quantum matrix algebra (by Example 2.6 and Theorem 3.5).
Note that the $4$-dimensional Sklyanin algebra and the multiparameter 
quantum matrix algebra satisfy the SSC of 
[{\bf Zh}, p. 398] (by [{\bf Zh}, 2.3(2)
and 2.4]) and hence satisfy the hypothesis $\Kdim M=\GKdim M$
in Theorem 3.5(3) (by [{\bf Zh}, 3.1]).

The Weyl algebras and the universal enveloping algebras 
of finite dimensional 
Lie algebras have standard filtrations, and 
their associated graded
rings are commutative polynomial rings. 
We prove the following theorem using similar ideas.

\proclaim{Theorem 0.3} Let $A_n$ be the $n$-th Weyl algebra with the 
standard filtration. Let $G$ be a finite group of filtered automorphisms
of $A_n$. If $|G|\neq 0$ in $k$, then $A_n^G$ is Auslander-Gorenstein
and $\GKdim$-Macaulay. Moreover if $\charr k=0$, then $A_n^G$ is 
Auslander-regular.
\endproclaim

A special case of Theorem 0.3 was proved by Levasseur, see [{\bf L1}, 3.2].
If $\charr k=0$, it is well-known that $A_n^G$ has finite global dimension
because $A_n$ is simple (see Remark 4.2). 
A quantum 
version of Theorem 0.3 is the following. We say that $q$ is {\it generic
with respect to $\{p_{ij}\}$} if $q^n$ ($n>0$) are not in the
the multiplicative subgroup generated by $\{p_{ij}\}$.

\proclaim{Theorem 0.4} Let $A_n(q,p_{ij})$ be the $n$-th quantum Weyl 
algebra defined in [{\bf GZ}, (2.3)] and $G$ a finite group of filtered 
automorphisms of $A_n(q,p_{ij})$. Suppose that $q$ is generic with 
respect to $\{p_{ij}\}$. If $|G|\neq 0$ in $k$, then $A_n(q,p_{ij})^G$ 
is Auslander-Gorenstein and $\GKdim$-Macaulay. 
\endproclaim

\proclaim{Theorem 0.5} Let $k$ be an algebraically closed field of
characteristic zero, $L$ a finite dimensional Lie algebra over $k$,
and $U(L)$ the universal enveloping algebra of $L$.
Let $G$ be a finite group of filtered automorphisms of $U(L)$. 
Then $U(L)^G$ is Auslander-Gorenstein and 
$\GKdim$-Macaulay if either of the  two following conditions holds:

(1) elements of $G$ act as inner automorphisms of the Lie algebra $L$;

(2) $L$ is one of the following simple Lie algebras:

$\quad$ (a) $A_n$, for $n=4m$ or $n=4m+1$;

$\quad$ (b) $B_n$;

$\quad$ (c) $C_n$;

$\quad$ (d) $E_6$, $E_7$, $E_8$, $F_4$, $G_2$.
\endproclaim

Theorem 0.5(1) was essentially known by Kraft and Small 
[{\bf KS}, Proposition 4(5)] (see Remark 6.7).

\msno
{\bf Definition 0.6.} Let $A$ be a noetherian ring. 

\msno
(1) The {\it grade} of an $A$-module $M$ is defined to be 
$$j(M)=\min\{i\;|\; \Ext^i_A(M,A)\neq 0\}$$
or $\infty$ if no such $i$ exists.

\msno
(2) We say $A$ is {\it
Auslander-Gorenstein} if the following two conditions hold:

(a) (Gorenstein property) $A$ has finite left and right injective dimension;

(b) (Auslander property) for every noetherian $A$-module $M$ and for all $i\geq 0$ and
all submodule $N\subset \Ext^i_A(M,A)$, $j(N)\geq i$.

We say $A$ is {\it Auslander-regular} if $A$ is Auslander-Gorenstein and
$A$ has finite global dimension.

\msno
(3) We say $A$ is {\it $\GKdim$-Macaulay} (where $\GKdim$ denotes the 
Gelfand-Kirillov dimension) if 
$$\GKdim\; M+ j(M)=\GKdim A<\infty$$
for every noetherian $A$-module $M$.

\medskip

Note that $\GKdim$-Macaulay is called Cohen-Macaulay 
by some authors.
Auslander-Gorenstein and $\GKdim$-Macaulay properties are considered as 
nice homological properties. Such rings are studied
by several researchers since 1980's (see [{\bf ASZ}], [{\bf B1}], [{\bf B2}], [{\bf BE}], [{\bf Ek}], [{\bf L1}], [{\bf L2}]
and so on). The skew polynomial rings,
the Weyl algebras and the universal enveloping 
algebras are Auslander--regular and $\GKdim$-Macaulay. 
Using the above
results one can construct more examples of Auslander-Gorenstein and 
$\GKdim$-Macaulay rings.

\bsno

\centerline{\bf 1. Preliminary}

\bsno

An old and basic tool
for filtered rings is to use their Rees rings
and
the associated graded rings. In this section we recall some  
definitions and properties about filtered rings and graded rings.
 
A {\it filtered 
algebra} is an algebra $A$ with a filtration $F=\{F_i\; |\; i\geq 0\}$ 
satisfying the conditions: (F0) $F_0=k$,\;  (F1) $F_i\subset F_j$ if 
$i\geq j$, \; (F2) $F_iF_j\subset F_{i+j}$, \; and (F3) $A=\bigcup_i F_i$.
Given a filtered algebra $A$, the Rees ring is defined to be 
$$\Rees A=\bigoplus_i F_n t^n.$$
The elements $t(=1t)$ and $t-1$ are central elements of $\Rees A$ and
$$\Rees A/(t-1)=A,\quad {\text{and}}\quad \Rees A/(t)=\G A.$$
Both $\Rees A$ and $\G A$ are {\it connected graded}  in the sense that the
degree zero part is $k$. By [{\bf L2}, 3.5], $\Rees A$ is noetherian if and only
if $\G A$ is noetherian. We call a filtered ring $A$ {\it filtered noetherian}
if $\G A$ is noetherian.

We will adopt the basic notations about graded rings
and graded modules from [{\bf JiZ}], [{\bf JoZ}] and [{\bf YZ}]. For example,
the opposite ring of $A$ is denoted by $A^{op}$.

Let $A$ be a connected graded ring. The maximal graded ideal $A_{\geq 1}$
is $A$ is denoted by $\fm$ and the trivial module $A/\fm$ is denoted by
$k$. The graded Ext-group is denoted by $\uExt$.

\msno
{\bf Definition 1.1.} Let $A$ be a noetherian connected graded ring.

\msno
(1) $A$ is called 
{\it AS-Gorenstein} (here AS stands for Artin-Schelter) if 

(a) $A$ is Gorenstein with injective dimension $d$, and

(b) there is an integer $l$ such that
$$\uExt^i_A(k,A)=\uExt^i_{A^{op}}(k,A)=
\left\{\matrix 0& {\text{for}} \; i\neq d\\
k(l) &{\text{for}} \; i=d
\endmatrix\right.$$
where $k(l)$ is the $l$-th degree shift of the trivial module $k$.

\msno
(2) $A$ is called 
{\it AS-regular} if it is AS-Gorenstein and has finite global dimension.

\medskip

Condition 1.1(1b) follows from 1.1(1a) if $A$ has enough normal elements
[{\bf Zh}, 0.2]. 
For example, the commutative polynomial rings are AS-regular. 
By [{\bf L2}, 6.3], a noetherian, connected graded, Auslander-Gorenstein ring
is AS-Gorenstein. 

Given a connected graded ring $A$ and any graded right $A$-module $M$, we
define the local cohomology of $M$ to be
$$H_{\fm}^i(M)=\varinjlim_n\uExt^i_A(A/A_{\geq n},M).$$
If $A$ is AS-Gorenstein, then $H_{\fm}^i(A)$ is 0 if $i\neq d$ and is
equal to $A'(l)$ if $i=d$. Here $(-)'$ is the graded vector space dual
of $-$.

Auslander property is also defined for non-Gorenstein ring if the ring
has a (rigid or balanced) dualizing complex [{\bf YZ}, 2.1]. In the case 
when $\G A$ is AS-Gorenstein, $A$ has Auslander property [Definition 0.6(2b)] 
is equivalent to that $A$ has an Auslander,
rigid dualizing complex [{\bf YZ}, 2.1 and 3.1]. We collect some results
below.

\proclaim{Theorem 1.2} (1) (Rees' lemma) If $A$ is a noetherian connected 
graded ring with a regular normal element $x$ of positive degree, then $A$ is 
AS-Gorenstein if and only if $A/(x)$ is.

(2) If $A$ is as in (1) $x$ is of degree 1, then $A$ is AS-regular 
if and only if $A/(x)$ is.

(3) If $A$ is a noetherian filtered algebra, then $\Rees A$ is AS-regular 
(respectively AS-Gorenstein) if and only if $\G A$ is.

(4) [{\bf L2}, 3.6] Let $A$ be a filtered noetherian algebra. Then 
$\Rees A$ is Auslander-Gorenstein (or Auslander-regular) if and only 
if $\G A$ is.

(5) [{\bf Ek}, 0.1; {\bf L2}, 3.6, 5.8 and 5.10] Let $A$ be a filtered 
noetherian algebra. If $\G A$ is Auslander-Gorenstein 
(respectively, and $\GKdim$-Macaulay) then so is $A$.
\endproclaim

We make the following definition. 

\msno
{\bf Definition 1.3.} Suppose $P$ is a property such that $\Rees A$ has $P$
if and only $\G A$ has $P$. We say a filtered algebra $A$ has {\it filtered
$P$} if $\Rees A$ (or $\G A$) has $P$.

\newpage
\centerline{\bf 2. Homological Determinant}

\bsno

First we recall some definitions in the graded case from [{\bf JiZ}] and [{\bf JoZ}]. 
Let $g\in \GrAut(A)$, the group of graded algebra 
automorphisms of $A$.
For any $A$-module $M$ we define 
the $g$-twisted module $M^g$ such that $M^g=M$ as a vector space and the
action is
$$m\ast a=mg(a)$$
for all $m\in M$ and $a\in A$. Let $f$ be a
$k$-linear homomorphism from $A$-module $M$ to $A$-module $N$. We
say $f$ is $g$-linear if
$$f(ma)=f(m)g(a)$$
for all $m\in M$ and $a\in A$. Then $f$ is a $g$-linear map if and only
if $f$ is an $A$-module homomorphism from $M$ to $N^g$.
Therefore $g$-linear maps can be extended to injective resolutions
(or projective resolutions). Moreover we can apply local cohomology
functor $H_{\fm}^\ast$
to $g$-linear maps. Now let $A$ be a graded AS-Gorenstein
ring with injective dimension $d$.  
By [{\bf JoZ}, 2.2 and 2.3], 
$g:A\to A$ induces a $g$-linear map $H_{\fm}^d(g): A'(l)\to
A'(l)$ where $l$ is the integer in Definition 1.1(1b) and
$'$ is the graded vector space dual, and $H_{\fm}^d(g)=c(g^{-1})'$. 
The constant $c^{-1}$ is called
the {\it homological determinant} of $g$
and we denote $\hdet_A\; g=c^{-1}$. By [{\bf JoZ}, 2.5],
$\hdet_A$ defines a group homomorphism $\GrAut(A)\to k^{\times}$.

The {\it trace} of $g$ on $A$ is defined to be
$$Tr_A(g,t)=\sum_{n\geq 0} {\operatorname{tr}}(g|_{A_n}) t^n.$$
See [{\bf JiZ}] for more details. The Hilbert series of $A$ is the trace of
the identity map, namely
$$H_A(t)=Tr_A(id, t)=\sum_{n\geq 0} \dim A_n\; t^n.$$
By [{\bf JoZ}, 2.6], if $A$ is AS-Gorenstein and $ g$ is 
rational in the sense of of [{\bf JoZ}, 1.3] then 
$$Tr_A(g,t)=(-1)^d(\hdet_A g)^{-1} t^{-l}+{\text{ lower
terms }},$$
when we expand the trace function  as a Laurent series
in $t^{-1}$. Here $d$ and $l$ are given in 
Definition 1.1(1b). Rationality of $g$ is automatic for 
AS-regular algebras [{\bf JoZ}, 4.2].

\proclaim{Lemma 2.1} (1) If $A$ is the commutative polynomial ring $k[V]$
and $g$ is a graded automorphism of $k[V]$, then $\hdet_A\; g=\det\; g|_V$
where $\det$ is the usual determinant of a $k$-linear map.

(2) If $A$ is the exterior algebra $\Lambda(V)$ and $g$ is a graded
automorphism of $A$, then $\hdet_A \; g=(\det\; g|_V)^{-1}$.
\endproclaim

\demo{Proof} (1) It follows by [{\bf Ben}, 2.5.1] that 
$Tr_A(g,t)=(\det \; (1-g|_Vt))^{-1}$. (Note that the action of $g$ in 
[{\bf Ben}] is defined via $g^{-1}$, see [{\bf Ben}, p. 1], so we change $g^{-1}$ in 
[{\bf Ben}, 2.5.1] to $g$.) Expand $(\det\; (1-gt)|_V)^{-1}$ as a Laurant series
in $t^{-1}$ we have
$$(\det\; (1-gt)|_V)^{-1}=(-1)^n
(\det\; g|_V)^{-1}t^{-n}+{\text{lower terms}}.$$
Therefore the result follows by [{\bf JoZ}, 2.6].

(2) It is easy to see that $Tr_A(g,t)=\det\; (1+g|_Vt)=(\det\; g|_V)t^n+
{\text{lower terms}}$ [{\bf Ben}, 5.2.1]. Therefore the result follows by
[{\bf JoZ}, 2.6].
\hfill $\square$
\enddemo

Lemma 2.1(2) shows that the homological 
determinant of $g$ may not
equal to the determinant of $g|_V$ (also see 
Examples 2.8 and 2.9). Lemma 2.1(1) shows that 
the homological determinant of $g$ is equal to the determinant of $g|_V$
when $A$ is the commutative polynomial ring.
The next example shows that Lemma 2.1(1) also holds for the $4$-dimensional
Sklyanin algebra. 

\medskip
\noindent
{\bf Example 2.2.} Let $k$ be the field of complex numbers $\Bbb C$.
Then the automorphisms of the $4$-dimensional Sklyanin algebra $S$ are
classified in [{\bf SS}, Section 2] and the generators of the automorphism
group are also listed there. 
So one can check $\det \; g|_V$ easily
where $V=S_1$ is the degree 1 part of $S$. By
[{\bf JoZ}, 2.6], the homological determinant of $g$ can be computed
by the trace of $g$. The trace of 
generators of the automorphism group is listed in [{\bf JiZ}, Section 4].
Use these facts
one can check that $\hdet_S \; g=\det \; g|_V$
for all $g\in \GrAut(S)$. 

\bigskip

We now extend the definition of homological
determinant to the filtered case.
Let $\FilAut(A)$ be the group of automorphisms of $A$ preserving the given
filtration $F$. Hence for any $g  \in \FilAut(A)$ one can extend $g  $
to a graded automorphism of $\Rees A$ by sending $t$ to $t$. Also $g  $
induces a graded automorphism of $\G A$. For simplicity we still use
$g  $ for both these graded automorphisms. Next lemma is clear.

\proclaim{Lemma 2.3} Let $B$ be a connected graded ring. Then $B=\Rees A$ 
for some filtered algebra $A$ if and only if there is a regular central 
element in $B_1$. 
\endproclaim 

The next proposition is useful for computing 
the homological determinant when $A$  has normal elements.

\proclaim{Proposition 2.4} 
Let $B$ be a noetherian graded AS-Gorenstein ring
and let $x$ be a regular normal element of positive degree. Suppose
$g  $ is a graded automorphism of $B$ such that $g  (x)=\lambda x$.
Then $\hdet_B\;  g  =\lambda \hdet_A\; g  $ where $A=B/(x)$.
\endproclaim

\demo{Proof} Suppose $d$ is the injective dimension of $B$. Then $d-1$ is
the injective dimension of $A$. Let $l_x: B(-s)\to B$ be the left
multiplication by $x$ where $s=\deg x$. Applying $H^\ast_{\fm}$ to the 
following commutative diagram
$$\CD
0@>>>    B(-s)       @>l_x>>   B @>>>   A @>>> 0\\
@.     @VV\lambda\; gV     @VVgV @VVgV  @.\\
0@>>>    B(-s)       @>l_x>>   B @>>>   A @>>> 0
\endCD
$$
we obtain that 
$$\CD
0@>>> H^{d-1}_{\fm}(A) @>>> H^{d}_{\fm}(B)(-s) @>>> H^{d}_{\fm}(B) @>>> 0\\
@.     @VV{{(\hdet_A g)^{-1}}(g^{-1})'}V @VV{{\lambda(\hdet_B g)^{-1}}(g^{-1})'}V @VVV  @.\\
0@>>> H^{d-1}_{\fm}(A) @>>> H^{d}_{\fm}(B)(-s) @>>> H^{d}_{\fm}(B) @>>> 0.
\endCD
$$
The result follows by comparing the first two vertical 
maps in the highest nonzero homogeneous component. 
\hfill $\square$
\enddemo

Let $A$ be a filtered noetherian AS-Gorenstein ring. 
If $g$ is a filtered automorphism of $A$, then $g(t)=t$
for $t\in \Rees A$ 
and by Proposition 2.4,
$\hdet_{\G A}g=\hdet_{\Rees A}g$. We define the 
{\it homological determinant} of $g$ on $A$ to be
$$\hdet_A g=\hdet_{\G A} g.$$

In the rest of this section we compute some homological
determinants.

\proclaim{Lemma  2.5} Let $A$ be a filtered algebra 
with $\G A=k_{p_{ij}}[x_1,\cdots,x_n]$, where $p_{ij}\neq p_{kl}$ for all
$(i,j)\neq (k,l)$. Let $g$ be a filtered automorphism of $A$.
Then $\hdet_A g=\det\; g|_{V}$ where $V=\bigoplus_{i=1}^n kx_i$. 
\endproclaim

\demo{Proof} First we assume that $A$ is $\Bbb N$-graded and $\deg x_i=1$
for all $i$.

Under the hypothesis on $p_{ij}$ it is routine to check that if 
$\sum_{i=1}^n a_ix_i$ is a normal element then it must be $a_ix_i$ for 
some $i$. Since $g(x_i)$ is normal, we can write $g(x_i)=b_ix_{\sigma(i)}$ 
for some permutation $\sigma\in S_n$ and $b_i\neq 0$. Applying $g$ to
the relations $x_ix_j=p_{ij}x_jx_i$ we see that
$p_{ij}=p_{\sigma(i),\sigma(j)}$. Therefore $\sigma$ is the identity 
by the hypothesis.

Now we consider the general case. By definition of $\hdet$ we can 
pass to the graded case. Let $g$ be a graded automorphism of $A=
k_{p_{ij}}[x_1,\cdots, x_n]$ with
$$\deg x_1\geq \cdots \geq \deg x_n\geq 1.$$
Let $w$ be the minimal index such that $\deg x_w=\deg x_n$. Then $g$ preserves 
$\oplus_{i=w}^n kx_i$. As an automorphism of $k_{p_{ij}}[x_w,
\cdots, x_n]$, $g(x_n)=b_nx_n$ for some
$b_n\neq 0$ by the previous argument.  By induction, $g(x_i)=b_ix_i+f_i(x_{i+1},\cdots, x_n)$
for some polynomial $f_i$.  
By Proposition 2.4, $\hdet_A g=b_n \hdet_{A/(x_n)} g$. Therefore 
$\hdet_A g=\prod_i b_i=\det g|_V$.  \hfill $\square$
\enddemo

\bsno
{\bf Example 2.6.} Let $M_{q,p_{ij}}(n)$ be the
multiparameter quantum matrix algebra. It is a connected graded
algebra generated by $\{x_{ij}\;|\; i,j=1,\cdots,n\}$ with
${{n^2}\choose 2}$ relations (see [{\bf AST}] for a list of relations),
where the structure constants $\{q, p_{ij}\}$ are generic.
Similar to the proof
of Lemma 2.5 with more careful computations we check that
$\hdet_{M_{q,p_{ij}}(n)} g=\det g|_V$ where $V=\bigoplus_{i,j=1}^n kx_{ij}$.
Details are left to interested readers.  

\proclaim{Proposition 2.7} Let $A$ and $B$ be two filtered AS-Gorenstein 
ring and assume that $A\otimes B$ is filtered noetherian. 
If $g$ (respectively $h$) is a filtered automorphism of $A$ (respectively
$B$), then $\hdet_{(A\otimes B)}(g\otimes h)=(\hdet_A g)(\hdet_B h)$.
\endproclaim

\demo{Proof} Let $\{F_i\;|\; i\geq 0\}$  (respectively $\{W_i\;|\; i\geq 0\}$)
be the filtration of $A$ (respectively $B$). The filtration of $A\otimes B$
is defined by $Z_n=\bigcup_{i+j=n} F_i\otimes W_j$. Then $\G (A\otimes B)
=(\G A)\otimes (\G B)$. By the definition of $\hdet$ in the filtered case, 
we may work with the associated graded rings, or equivalently we may 
assume $g$ and $h$ are graded
automorphisms of graded AS-Gorenstein rings $A$ and $B$ respectively. 
Since systems $\{(A\otimes B)_{\geq n}\}$ and $\{A_{\geq n}\otimes B+A\otimes
B_{\geq n}\}$ are cofinal, we have
$$H^{n}_{\fm_{A\otimes B}}(M_A\otimes N_B)=\bigoplus_{i+j=n} 
H^{i}_{\fm_{A}}(M_A)\otimes H^{j}_{\fm_{B}}(N_B).$$
In particular, if $E_A$ and $F_B$ are graded injective modules over
$A$ and $B$ respectively, then $E_A\otimes F_B$ is 
$H^{\ast}_{\fm_{A\otimes B}}$-acyclic. 
Let $E^\ast$ and $F^\ast$ be injective
resolutions of $A$ and $B$ respectively. We can compute $H_{\fm}^\ast(A
\otimes B)$ by the resolution $E^\ast\otimes F^\ast$. Since $g$ gives a 
$g$-linear map of $E$ and $h$ gives an $h$-linear of $F$, $g\otimes h$
gives a $g\otimes h$-linear map of $E^\ast\otimes F^\ast$.
Therefore
$$H^{d+e}_{\fm_{A\otimes B}}(g\otimes h)=H^{d}_{\fm_A}(g)\otimes H^{e}_{\fm_B}
(h),$$
where $d=\id A$ and $e=\id B$. By the definition of homological determinant
in the graded 
case, $\hdet_{A\otimes B} (g\otimes h)=(\hdet_A g)(\hdet_B h).$
\hfill $\square$
\enddemo

Finally we give two examples where $\hdet_A g$ is not equal to
$\det g|V$ where $V$ is the minimal generating space.

\bsno
{\bf Example 2.8.} Let $A$ be the commutative  graded
ring $k[x_1,x_2]/(x_1^3)$. Let $g$ be the graded
automorphism of $A$ sending $x_i$ to $-x_i$ for $i=1,2$.
Then $V=kx_1\oplus kx_2$ is the minimal generating space
and $\det g|V=1$. Also $g^2=id$. Since $g(x_1^3)=-x_1^3$,
by Proposition 2.4, $\hdet_A g=(-1)\hdet_{k[x_1,x_2]} g=-1$.
Hence $\hdet_A g$ is neither equal to $\det g|_V$ nor $(\det g|_V)^{-1}$.
It is also easy to check that $A^G$ is not Gorenstein
with $G=\{g,id\}$.

\bsno
{\bf Example 2.9.} Let $A$ the skew polynomial ring generated by $x, y, z$ 
with relations $xy=-yx, xz=zx, yz=zy$. Let $g$ be the graded automorphism
of $A$ such that $g(x)=-y, g(y)=-x, g(z)=-z$. Then $g^2=1$
and $\det g|_V=-1$. The Koszul dual is  $A^!=k\langle x, y, z\rangle/I$, 
where $I$ is generated by $x^2=y^2=z^2=0, xy=yx, yz=-zy, xz=-zx$. The
induced automorphism $g^{\tau}$ of $A^!$ is given by $x\to -y, y\to -x, 
z\to -z$. It follows from [{\bf JiZ}, 3.4] that
$$Tr_A(g, t)=\frac 1{Tr_{A^!}(q^{\tau}, -t)}=\frac 1{1+t+t^2+t^3}.
$$
Therefore $\hdet_A g=1$ [{\bf JoZ}, 2.6]. Let $G= \{1, g\}$. Then 
the Hilbert series of $A^G$ is
$$H_{A^G}(t)=\frac 12(Tr_A(1, t)+Tr_A(g, t))
\neq \pm t^m H_{A^G}(t^{-1})$$
for any $m\in {\Bbb Z}$. So $A^G$ is not Auslander-Gorenstein
[{\bf Ei}, Ex. 21.17(c)].
This gives another example of $\hdet g\neq \det g|_V$ and
also shows that Theorem 0.1 fails without the extra hypothesis
on $p_{ij}$.
\bsno

\centerline{\bf 3. Invariant subring of filtered rings}

\bsno

In this section we prove Theorems 0.1 and 0.2. The following lemma is clear.

\proclaim{Lemma 3.1} Let $A$ be a filtered ring and 
$G$ be a subgroup of $\FilAut(A)$. Then

(1) $\Rees A^G=(\Rees A)^G$.

(2) If $|G|\neq 0$ in $k$, then $\G A^G=(\G A)^G$.
\endproclaim

If $|G|=0$ in $k$, then Lemma 3.1(2) may fail.

\bsno
{\bf Example 3.2.} 
Let $k$ be a field of characteristic $2$ and let $g$ be a filtered 
automorphism of the polynomial ring $k[x]$ determined by
$g(x)=x+1$. Hence $g^2=id$. Let $G=\{id,g\}$. Consider $k[x]$ as a 
filtered ring and $G$ a subgroup of filtered automorphisms. Then 
the induced graded automorphism of $g$ is the identity on the $k[x]
=\G k[x]$. Hence $k[x^2]=\G A^G\subsetneq (\G A)^G=k[x]$. 

\medskip

We now recall results from [{\bf JoZ}] and [{\bf YZ}]. 

\proclaim{Theorem 3.3} Let $A$ be a noetherian connected graded ring
and $G$ a finite subgroup of $\GrAut(A)$ with $|G|\neq 0$
in $k$. Suppose that $\hdet_A g=1$ for all $g\in G$. 

(1) [{\bf JoZ}, 3.3] If $A$ is AS-Gorenstein, then so is $A^G$.

(2) [{\bf YZ}, 4.20] If $A$ is Auslander-Gorenstein, then so is
$A^G$.

(3) [{\bf YZ}, 5.13 and 5.14] If $A^G$ is AS-Gorenstein and $A^G$ either is
PI or has enough normal elements, then 
$A^G$ is Auslander-Gorenstein and $\GKdim$-Macaulay.
\endproclaim

A slightly more general version of Theorem 3.3(3) holds.
Let $\Kdim$ denote the Krull dimension. 

\proclaim{Lemma 3.4} Let $B\subset A$ be two noetherian connected
graded rings such that $A$ is finitely generated over $B$ on both sides
and $B$ is a direct summand of $A$ as graded $B$-bimodules. Suppose
$\GKdim M=\Kdim M<\infty$ for all finitely generated graded $A$-module
$M$. Then 
$\GKdim N=\Kdim N=\Kdim N\otimes_B A$
for all finitely generated graded $B$-module $N$.
\endproclaim

If $G$ is a finite subgroup of $\GrAut(A)$ with $|G|\neq 0$ in $k$,
then $B:=A^G\subset A$ is a direct summand of $A$. If moreover
$A$ is noetherian, then $A$ is finitely generated over $B$ on both sides
[{\bf Mo}, 2.1.1 and 1.12.1].

\demo{Proof} Since $B$ is a direct summand of $A$ as $B$-bimodule,
$N_B$ is a direct summand of $N\otimes_B A_B$. Hence
$\GKdim N\leq \GKdim N\otimes_B A_B$ and $\Kdim N\leq \Kdim N\otimes_B A_B$.

By [{\bf MR}, 8.3.14(iii)], $\GKdim N\geq \GKdim N\otimes_B A_A$ and
$\GKdim M_A\geq \GKdim M_B$ for any $A$-module $M$. Hence
$$\GKdim N\otimes_B A_A= \GKdim N\otimes_B A_B= \GKdim N.$$
It remains to show that $\Kdim N=\Kdim N\otimes_B A$. For every submodule
$L\subset N$, let $\phi$ be the map $L\otimes_B A\to N\otimes_B A$.
Then $L\to \phi(L\otimes_B A)$ is an order-preserving map from
the lattice of $B$-submodules of $N$ to the lattice of $A$-submodules
of $N\otimes_B A$. Since $A=B\oplus C$ for some $C$ as $B$-bimodule,
this map preserves proper containment. Therefore
$$ \Kdim N\leq \Kdim N\otimes_B A.$$
We now prove $\Kdim N\geq \Kdim N\otimes_B A(=\GKdim N)$ by induction
on $\Kdim N$. By noetherian induction we may assume $N$ is 
$\Kdim$-critical. Let $s=\Kdim N\otimes_B A=\GKdim N<\infty$. Pick an 
$A$-submodule $M\subset N\otimes_B A$ such that $(N\otimes_B A)/M$
is critical of Krull (and GK) dimension $s-1$ and let $\bar{N}$ be the image
of the composition
$$N\to N\otimes_B A\to (N\otimes_B A)/M.$$
Then
$$\GKdim \bar{N}\leq \GKdim \bar{N}A_B\leq \GKdim \bar{N}\otimes_B A_B
=\GKdim \bar{N}.$$
This implies that the first $``="$ of the following
$$\GKdim \bar{N}A_A\geq \GKdim \bar{N}A_B=\GKdim \bar{N}
=\GKdim \bar{N}\otimes_B A_A\geq \GKdim \bar{N}A_A.$$
Therefore
$$\GKdim \bar{N}= \GKdim \bar{N}A_A=
\GKdim (N\otimes_B A)/M=s-1.$$
Hence $\bar{N}$ is a proper quotient of $N$ and by induction hypothesis
$\Kdim \bar{N}\geq \GKdim \bar{N}=s-1$. Since $N$ is critical,
$\Kdim N\geq \Kdim \bar{N}+1\geq s$ as required. \hfill $\square$
\enddemo

By [{\bf Zh}, 2.3, 2.4 and 3.1], if $A$ has enough normal elements or if $A$ is the
$4$-dimensional Sklyanin algebra (or more generally if $A$ satisfies SSC)
then $$\Kdim M=\GKdim M$$ for all finitely generated graded $A$-module $M$.

\proclaim{Theorem 3.5} Let $A$ be a noetherian filtered ring and $G$ be a 
finite subgroup of $\FilAut(A)$. Suppose that $|G|\neq 0$ in $k$ and
that $\hdet_A g=1$ for all $g\in G$. 

(1) If $A$ is filtered AS-Gorenstein,
then so is $A^G$.

(2) If $A$ is filtered Auslander-Gorenstein, 
then so is $A^G$.

(3) Suppose $\GKdim M=\Kdim M$ for all finitely generated
graded $\G A$-module $M$.
If $A$ is filtered Auslander-Gorenstein 
and $\GKdim$-Macaulay, then so is $A^G$.
\endproclaim

\demo{Proof} (1) By definition we only need to show $\G A^G$ is 
AS-Gorenstein. 
Now the result follows by Theorem 3.3(1).

\noindent
(2) It suffices to show $\G A^G$ is Auslander-Gorenstein,
which follows by Theorem 3.3(2).

\noindent
(3) We can pass to the graded case and assume that $A$ is graded. 
By (2) $A^G$ is Auslander-Gorenstein. It remains to show
the $\GKdim$-Macaulay property. Let $\Cdim$ be the canonical dimension 
defined in [{\bf YZ}, 2.9]. Since $A$ is $\GKdim$-Macaulay, $\GKdim M=\Cdim
M$  for all finitely generated graded $A$-module $M$. Let $B$ denote
$A^G$. For every finitely generated graded $B$-module $N$, it follows 
from Lemma 3.4 that
$$\Kdim N=\GKdim N=\GKdim N\otimes_{B} A_{A}$$
By [{\bf YZ}, 4.14], $\Kdim N \leq \Cdim N$ and it is clear that
$$\Cdim N\leq \Cdim N\otimes_{B} A_{B}=\Cdim N\otimes_{B} A_{A}
= \GKdim N\otimes_{B} A_{A}$$
where the first $``="$ is [{\bf YZ}, 4.17(3)].
Combining these (in)equalities we obtain $ \Kdim N= \GKdim N=\Cdim N$.
Therefore $A^G$ is graded $\GKdim$-Macaulay. By [{\bf L2}, 5.8], $A^G$
is (ungraded) $\GKdim$-Macaulay as required.
\hfill $\square$
\enddemo

The following example is well-known, which shows that the hypothesis
$\hdet_A g=1$ can not be deleted from Theorem 3.5.

\msno
{\bf Example 3.6.} 
Let $A$ be the skew polynomial ring ${\Bbb C}_{p_{ij}}
[x,y,z]$ \, where $p_{ij}\neq 0\in \Bbb C$,
and $g$ a graded
automorphism sending $x\to -x$, $y\to -y$ and $z\to -z$. Then $g^2=id$
and $\det g=-1$. By Proposition 2.4, $\hdet g=-1$ and
it is easy to check that $Tr(g, t)=\frac 1{(1+t)^3}$. Let $G=\{id, g\}$.
Then the Hilbert series of $A^G$ is 
$$H_{A^G}(t)={1\over 2}({1\over (1-t)^3}+{1\over (1+t)^3})
={{1+3t^2}\over{(1-t^2)^3}}.$$
So $H_{A^G}(t^{-1})\neq \pm t^m H_{A^G}(t)$ for any $m\in {\Bbb Z}$.
Therefore $A^G$ is not Gorenstein by [{\bf Ei}, Exercise 21.17(c)].

\medskip

To apply Theorem 3.5, one needs to show $\hdet g=1$ for all $g\in G$.
To show $\hdet g=1$ it suffices to show the only group homomorphism
$G\to k^\times:=k-\{0\}$ is trivial. Here are some easy cases.

\proclaim{Lemma 3.7} [{\bf JoZ}, 3.4] 
The homomorphism $G\to k^\times$ is trivial if
one of the following holds:

(a) $G=[G,G]$;

(b) $|G|$ is odd and $k={\Bbb Q}$;

(c) $k={\Bbb Z}/2{\Bbb Z}$.
\endproclaim

One can easily get some corollaries by combining Lemma 3.7 and Theorem 3.5.
We now prove Theorems 0.1 and 0.2.

\demo{Proof of Theorem 0.1} By [{\bf Zh}, 2.3(2) and 3.1], $\GKdim M=\Kdim M$
holds for all finitely generated graded 
$k_{p_{ij}}[x_1,\cdots, x_n]$-module $M$. So the
result follows from Lemma 2.5 and Theorem 3.5(3).
\hfill $\square$
\enddemo

\demo{Proof of Theorem 0.2} Again by [{\bf Zh}, 2.3(2) and 3.1], 
$\GKdim M=\Kdim M$ holds for all
finitely generated graded $k_{p_{ij}}[x_1,\cdots, x_n]\otimes 
\Lambda(V)$-module $M$. 

Since $A$ is a filtered algebra, $A\otimes \Lambda(V)$
is a filtered algebra with the filtration induced by the filtration
of $A$ and $\G (A\otimes \Lambda(V))=(\G A)\otimes \Lambda(V)
=k_{p_{ij}}[x_1,\cdots, x_n]\otimes \Lambda(V)$. By [{\bf Zh}, 0.2]
and Theorem 1.2(5), $A\otimes \Lambda(V)$ is Auslander-Gorenstein and
$\GKdim$-Macaulay. Also $A\otimes \Lambda(V)$ and 
$(A\otimes \Lambda(V))^G$ are filtered noetherian. 

By Proposition 2.7 and Lemmas 2.5 and 2.1(2),
$$\hdet_{A\otimes \Lambda(V)} g
=\hdet_A g\cdot\hdet_{\Lambda(V)} g=\det g|_V\cdot (\det g|_V)^{-1}
=1$$ 
for all $g\in G$. Hence the result follows from Theorem 3.5(3).
\hfill $\square$
\enddemo

Let $\Lambda_{p_{ij}}(V)=k\langle x_1, \cdots, x_n\rangle/ 
(x_i^2, x_ix_j+p_{ij}^{-1}x_jx_i)$ be the quantized exterior algebra. 
By Theorem 3.5(3) and similar
arguments as in the proof of Theorem 0.2 we have 
the following result which contains the commutative case [{\bf Ben}, 5.3.2]
as a special case.

\proclaim{Theorem 3.8} 
Suppose $A$ is a filtered ring such that the associated
graded ring $\G A$ is isomorphic to the skew polynomial ring 
$k_{p_{ij}}[x_1,\cdots, x_n]$. Let $G$ be a finite
group of filtered automorphisms of $A$ with $|G|\neq 0$ in $k$. Then 
$(A\otimes \Lambda_{p_{ij}}(V))^G$ 
is Auslander-Gorenstein and $\GKdim$-Macaulay 
where $V=\bigoplus_{i=1}^n kx_i$.
\endproclaim

The following is an immediate consequence of Theorems
0.1-0.2 and [{\bf YZ}, 6.23].

\proclaim{Corollary 3.9} Let $A^G$ be as in
Theorem 0.1 and $(A\otimes \Lambda(V))^G$ be as in 
Theorem 0.2. Then they have quasi-Frobenius rings of fractions.
\endproclaim

Another method to check whether $A^G$ is Gorenstein is to use a noncommutative
version of Stanley's theorem [{\bf JoZ}, Sect. 6]. For example if $\G A$ 
(or $\Rees A$) is Auslander-regular and $\charr k=0$, then $A^G$ is 
Auslander-Gorenstein if the Hilbert series of $\G A^G$ satisfies 
$$H_{\G A^G}(t^{-1})=\pm t^m H_{\G A^G}(t)$$
for some $m\in {\Bbb Z}$ (see [{\bf JoZ}, 6.2 and 6.4]).

\bsno

\centerline{\bf 4. Invariant subrings of the Weyl algebras}

\bsno

Let $A_n$ be the $n$-th Weyl algebra and $G$ a finite group of automorphisms
of $A_n$. If $\charr k=0$, then it is well-known that $A_n^G$ has 
global dimension $n$ for any $G$. We start with 
some basic facts below.

\proclaim{Theorem 4.1} Let $\charr k=0$ and $G$ a finite group of
automorphisms of the $n$-th Weyl algebra $A_{n}$. Then $A_n^G$
is a noetherian simple domain of global dimension and Krull dimension $n$
and Gelfand-Kirillov dimension $2n$.
\endproclaim

\demo{Proof} By [{\bf Mo}, 1.12.1], $A_n^G$ is 
noetherian. It is clear that $A_n^G$ is a domain. 

Recall that an automorphism $g$ is inner if there is a
unit $u\in A$ such that $g(x)=uxu^{-1}$. Since units of $A_{n}$
are elements of $k^{\times}$, every automorphism of $A_n$ is outer
(namely, not inner).
By [{\bf Mo}, 2.6] and [{\bf Mo}, 2.4], $A_n^G$ is simple and $A_n$ is a finitely
generated projective $A_n^G$-module. Hence the standard spectral 
sequence
$$E^{p,q}_2:=\Ext^p_{A_n}(\Tor^{A_n^G}_q(N,A_n),M)\Rightarrow 
\Ext^{p+q}_{A_n^G}(N, M)$$
collapses to the isomorphisms
$$\Ext^p_{A_n}(N\otimes_{A_n^G}A_n,M)=\Ext^{p}_{A_n^G}(N, M)$$
for all right $A_n$-module $M$ and right $A_n^G$-module $N$.
Therefore $\Ext^{p}_{A_n^G}(N, M)=0$ for all $p>n$.
For any right $A_n^G$-module $M_0$, let $M=M_0\otimes_{A_n^G}A_n$.
Since $A_n^G$ is an $(A_n^G,A_n^G)$-bimodule direct summand of $A_n$
[{\bf Mo}, 2.1.1],
$M_0$ is a direct summand of $A_n^G$-module $M$. Thus
$\Ext^{p}_{A_n^G}(N, M_0)=0$ for all right $A_n^G$-modules $N$ and $M_0$.
This implies that $\gldim A_n^G\leq n$. Pick a simple $A_n$-module
$N$ with $\GKdim n$, then $\GKdim_{A_n} N\otimes_{A_n^G} A_n=n$ 
[{\bf MR}, 8.3.14(iii)], thus $\Ext^n_{A_n}(N\otimes_{A_n^G} A_n, A_n)\neq 0$. 
Consequently, $\Ext^n_{A_n^G}(N, A_n)\neq 0$ and $\gldim A_n^G\geq n$.
So $\gldim A_n^G=n$. 

For GK-dimension we have
$$2n=\GKdim A_n=\GKdim (A_n^G\otimes_{A_n^G} A_n)\leq \GKdim A_n^G
\leq \GKdim A_n$$
where the first $\leq$ is [{\bf MR}, 8.3.14(iii)]. Therefore $\GKdim A_n^G=2n$.

Since $A_n$ is projective over $A_n^G$ and contains $A_n^G$ as a
direct summand, $A_n$ is faithfully  flat over $A_n^G$. Thus
$\Kdim A_n^G\leq\Kdim A_n$. For any group action one can form a skew group
ring $A_n\ast G$ which is free (hence faithfully flat) over $A_n$. 
Thus $\Kdim A_n\leq \Kdim A_n\ast G$. By [{\bf Mo}, 2.5], $A_n^G$
and $A_n\ast G$ are Morita equivalent. Therefore
$\Kdim A_n^G=\Kdim A_n\ast G=\Kdim A_n=n$. 
\hfill $\square$
\enddemo

\noindent
{\bf Remark 4.2.} The above proof shows the following known result:
If $A$ is noetherian and simple with finite global dimension  and $G$ is 
a finite group of (outer) automorphisms of $A$ with $|G|\neq 0$ in $k$, then 

(a) $A^G$ is noetherian and simple;

(b) $\gldim A^G=\gldim A$;

(c) $\Kdim A^G=\Kdim A$;

(d) $\GKdim A^G=\GKdim A$.

\medskip

It is not clear from the above proof that $A_n^G$ satisfies Auslander and
and $\GKdim$-Macaulay properties (which we believe to be true). 
Next we are going
to show for certain automorphism group $G$ (namely for filtered 
automorphism group) $A_n^G$ satisfies the Auslander  and 
$\GKdim$-Macaulay properties without any restriction on $\charr k$.

Let $A_n$ be generated by $x_1,\cdots, x_n$
and $y_1,\cdots, y_n$ subject to the relations
$$[x_i,x_j]=[y_i,y_j]=0, \quad [y_i,x_j]=\delta_{ij}, \quad {\text{for all}}
\quad i,j.$$
This is a filtered algebra with the standard filtration determined by
$F_1=k+ \sum_i kx_i+\sum_j ky_j$. The associated graded ring
is the commutative polynomial ring $k[x_1,\cdots, x_n,y_1,\cdots, y_n]$.

Let $V_{n}$ be the vector space generated by these $x_1,\cdots, x_n$ 
and $y_1,\cdots, y_n$. We say a linear map $\sigma:V_{n}\to V_{n}$ 
is a $[-,-]$-map if 
$$[\sigma(x_i),\sigma(x_j)]=[\sigma(y_i),\sigma(y_j)]=0, \quad 
[\sigma(y_i),\sigma(x_j)]=\delta_{ij}, \quad {\text{for all}}
\quad i,j.$$
We call another basis $\{x'_1,\cdots, x'_n, y'_1,\cdots y'_n\}$ of $V_{n}$
an $[-,-]$-basis if the map $\sigma:x_i\to x'_i, y_i\to y'_i$
is a $[-,-]$-map.

\proclaim{Lemma 4.3} If $g$ is a filtered endomorphism of $A_n$, then 
there is a $[-,-]$-map
$\sigma: V_{n}\to V_{n}$ and a linear map $\epsilon: V_{n}\to k$
such that $g|_{V_{n}}=\sigma+\epsilon$. Further if $\sigma$ is an 
automorphism, then $g$ is a filtered automorphism of $A_n$ and
$\hdet_{A_n} g =\det \sigma|_{V_{n}}$.
\endproclaim 

\demo{Proof} Since $g$ maps $F_1=k\oplus V_{n}$ to $F_1$, the restriction
$g|_{V_{n}}$ decomposes into two parts $\sigma$ and $\epsilon$.
For every $a,b\in V_{n}$, $[a,b]\in k$ and hence 
$$[a,b]=g([a,b])=[g(a),g(b)]=[\sigma(a)+\epsilon(a), \sigma(b)+\epsilon(b)]
=[\sigma(a), \sigma(b)]$$
because $\epsilon(a)$ commutes with all elements. Therefore
$\sigma$ is a $[-,-]$-map. If $\sigma$ is an automorphism then
$g$ is an automorphism of $F_1$ and hence a filtered automorphism of
$A_n$. By Lemma 2.1(1),
$$\hdet_{A_n} g=
\hdet_{k[V_{n}]} g=\det g|_{V_{n}}=\det \sigma|_{V_{n}}.$$
\hfill $\square$
\enddemo

\proclaim{Lemma 4.4} If $\sigma: V_{n}\to V_{n}$ is a $[-,-]$-map,
then $\det \sigma=1$. As a consequence, $\sigma$ is an automorphism.
\endproclaim

\demo{Proof} Replacing $k$ by its algebraic closure will not change the
determinant. So we may assume that $k$ is algebraically closed. Let $x$
be an eigenvector of $\sigma$. So $\sigma(x)=s x$ for some $s\in k$. Write
$x=\sum a_i x_i+\sum b_i y_i$. We may assume some $a_i\neq 0$. We will 
make several base changes such that $x=x_n$. Apparently we
will require that all base changes preserve the bracket relations.
By changing $x_i$ and then changing $y_i$ properly (to keep the $[-,-]$),
we have another $[-,-]$-basis so that $x=x_n+\sum b_i y_i$. By changing
$\{y_1,\cdots, y_{n-1}\}$ we may assume $x=x_n+b_{n-1}+b_n y_n$. Changing
basis within $\{x_n, y_n\}$, we have $x=x_n+b_{n-1}y_{n-1}$. Exchanging
$x_{n-1}$ and $-y_{n-1}$, we have $x=x_n+b_{n-1}x_{n-1}$, Finally we
change $x_i$ so that $x=x_n$. Since $\sigma$ is a $[-,-]$-map,
$$1=[y_n,x_n]=[\sigma(y_n),\sigma(x_n)]=[\sigma(y_n),sx_n].$$
Thus $s\neq 0$ and $\sigma(y_n)=s^{-1}y_{n}+\sum_{i} a_ix_i+\sum_{i<n}b_iy_i$.
Let $T=\bigoplus_{i=1}^n kx_i+\bigoplus_{i<n} ky_i$. Then $T=
\{x\in V_{n}\;|\; [x_n,x]=0\}$. Since $\sigma$ is a $[-,-]$-map, $T$ is
$\sigma$-invariant, namely, $\sigma(T)\subset T$. Write $T=V_{(n-1)}\oplus
ky_n$. We may decompose the restriction $\sigma|_{V_{(n-1)}}$ into two 
linear maps $\sigma|_{V_{(n-1)}}=\theta+\eta$ where $\eta: V_{(n-1)}
\to ky_n$ is a linear map and $\theta: V_{(n-1)}\to 
V_{(n-1)}$ is a $[-,-]$-map because $y_n$ commutes with $V_{(n-1)}$.
By induction $\det \theta=1$ and clearly 
$$\det \sigma=s\times s^{-1}\times \det \theta=1.$$
As a consequence $\sigma$ is invertible.
\hfill $\square$
\enddemo

By Lemmas 4.3 and 4.4 the following is clear.

\proclaim{Corollary 4.5} Every filtered endomorphism of the Weyl algebra
$A_n$ is an automorphism.
\endproclaim 

\proclaim{Theorem 4.6} Let $G$ be a finite group of filtered automorphisms
of the Weyl algebra $A_n$. If $|G|\neq 0$ in $k$, then $A_n^G$ is filtered 
Auslander-Gorenstein and $\GKdim$-Macaulay. 
\endproclaim

\demo{Proof} This is a consequence of Theorem 0.1 and Lemmas 4.3 and 4.4.
\hfill $\square$
\enddemo

Theorem 0.3 is an immediate consequence of Theorems 4.1 and 4.6.

%
%
%
%

For the first Weyl algebra $A_1$ over the complex numbers $\Bbb C$, 
it is known that every finite subgroup of $\Aut(A_1)$ is conjugate to 
a subgroup of $SL(2,{\Bbb C})\subset \FilAut(A_1)$ (see [{\bf AHV}]).
A complete list of finite subgroup of $SL(2,{\Bbb C})$ is also 
listed in 
[{\bf AHV}]. Therefore we have the following. 

\proclaim{Corollary 4.7} Let $A_1$ be the first Weyl algebra over $\Bbb C$
and $G$ a finite group of automorphisms of $A_1$. Then $A_1$ is
Auslander-regular and $\GKdim$-Macaulay.
\endproclaim

\bsno

\centerline{\bf 5. Invariant subrings of quantum Weyl algebras}

\bsno

In this section we look at  a family of quantum Weyl algebras
studied in [{\bf GZ}] and prove similar statements as in Section 4. 

Let $\{q\}\cup \{p_{ij}\;|\; i,j=1,\cdots, n\}$ be a set of nonzero 
elements in $k$ with $p_{ii}=1$ and $p_{ij}=p_{ji}^{-1}$. The quantum
Weyl algebra $A_n(q,p_{ij})$ is generated by the elements
$x_1,\cdots, x_n$, $y_1,\cdots, y_n$ subject to the relations
$$\align 
x_ix_j& =p_{ij}qx_jx_i, \qquad \quad {\text{for all}} \quad \; i<j\\
y_iy_j& =p_{ij}q^{-1}y_jy_i, \qquad {\text{for all}} \quad \; i<j\qquad\qquad \qquad\qquad\qquad(E1)\\
y_ix_j& =p^{-1}_{ij}q x_jy_i, \qquad \;\; {\text{for all}} \quad \; i\neq j\\
y_ix_i&=1+q^2 x_iy_i+(q^2-1)\sum_{j>i}x_jy_j, \qquad {\text{for all}} \quad
\; i.
\endalign
$$
By [{\bf GZ}, 3.11(1)], $A_n(q,p_{ij})$ is a noetherian, Auslander-regular and 
$\GKdim$-Macaulay domain of GK-dimension $2n$. It is easy to check that
the associated graded ring has enough normal elements. Hence it follows 
from [{\bf Zh}, 2.3(2) and 3.1] that the hypothesis $\GKdim M=\Kdim M$ 
in Theorem 3.5(3) is satisfied. If $q^2=1$ and $\charr k=0$, 
then it is simple and has global dimension and Krull dimension $n$ 
[{\bf FKK}]. If $q^2\neq 1$ or $\charr k\neq 0$, then 
$A$ is not simple and has global dimension and Krull dimension $2n$ 
[{\bf GZ}, 3.11(3)]. If $q^n\neq 1$ for all $n>0$, then it is primitive [{\bf GZ}, 3.2].

\proclaim{Lemma 5.1} Suppose that $q^4\neq 1$, 
$q p_{ij}\neq 1$, $q^3p_{ij}\neq 1$ for all $i\neq j$.
Then 

(1) Every filtered automorphism $g$ of $A_n(q,p_{ij})$ 
has the form
$$g(x_i)=\alpha_ix_i+\sum_{j>i} a_{ij}x_j+
\sum_{j>i} b_{ij}y_j$$ 
$$g(y_i)=\alpha_i^{-1}y_i+\sum_{j>i} c_{ij}x_j+
\sum_{j>i} d_{ij}y_j$$ 
for $\alpha_i\neq 0, a_{ij}, b_{ij}, c_{ij}, d_{ij}\in k$.

(2) $\det g|_{V_n}=1$ where $V_n=\bigoplus_{i=1}^n
 kx_i\oplus \bigoplus_{i=1}^n ky_i$. 
\endproclaim

\demo{Proof} (2) It follows from (1) immediately.

\noindent
(1) By [{\bf GZ}, 1.5] and the relations above, the ordered monomials
in $x_1,\cdots, x_n, y_1,\cdots, y_n$ are basis of $A_n(q, p_{ij})$.
Write 
$$\align
g(x_n)&=\sum_{j}a_{j}x_j+\sum_j b_{j}y_j+t\\
g(y_n)&=\sum_{j}c_{j}x_j+\sum_j d_{j}y_j+s.
\endalign
$$ 
Use the relations listed in (E1) and
$$g(y_n)g(x_n)=1+q^2g(x_n)g(y_n) \qquad\qquad (E2)$$
we obtain $a_ic_i=q^2 a_ic_i$ by comparing the 
coefficients in $x_i^2$ terms. If $a_{i_0}\neq 0$ for some 
$i_0$ then $c_{i_0}=0$ since $q^2\neq 1$. By comparing 
the coefficients of $x_{i_0}x_j$ terms for $j\neq i_0$, we obtain
$c_ja_{i_0}p_{ji_0}q^{-1}=a_{i_0}c_j q^2$. Thus $c_j=0$ since
$q^3p_{i_0j}\neq 1$. Since $g$ is an algebra automorphism
$d_j\neq 0$ for some $j$. A similar argument by exchanging
$x's$ and $y's$ we see that $b_j=0$ for all $j$. Thus (E2)
implies that $s=t=0$. If there are $a_i\neq 0$ and 
$d_j\neq 0$ for $j\neq i$ then $a_id_j=d_ja_i p_{ij}q$
by comparing the $x_iy_j$ terms. This
contradicts with $p_{ij}q\neq 1$. Thus $a_j=d_j=0$
for all $j\neq i_0$. Now if we compare (E2) with (E1) we 
see that $g(x_n)=a_n x_n$ and $g(y_n)=a_{n}^{-1} y_n$ 
and write $\alpha_n=a_n$. In the other case when $a_i=0$ for
$i$, we can show that this is impossible similar to the above
argument (the condition $q^4\neq 1$ will be used in this case). Thus
$$g(x_n)=\alpha_n x_n, \qquad g(y_n)=\alpha_n^{-1} y_n.$$

Now we decompose the vector space $V_n$ into $V_{n-1}\bigoplus
(kx_n\oplus ky_n)$ and decompose $g|_{V_n}$ into
$$g|_{V_{n-1}}=g'|_{(V_{n-1}\to V_{n-1}\oplus k)}+ 
\epsilon|_{(V_{n-1}\to kx_n\oplus ky_n)}.$$
It is not hard to see that $g'$ can be extended to 
an algebra automorphism of $A_{n-1}(q,p_{ij})$. By induction
$g'$ has the required form and the statement follows 
by induction. \hfill $\square$
\enddemo

\demo{Proof of Theorem 0.4} If $q$ is generic, then 
the hypothesis of Lemma 5.1 holds. The result follows from
Theorem 3.5(3) and Lemma 5.1. \hfill $\square$
\enddemo

\noindent
{\bf Remark 5.2}. Theorem 0.4 should also hold even 
if $q$ is not generic and holds for other nice quantum
Weyl algebras studied in [{\bf GZ}]. 

\bsno

\centerline{\bf 6. Invariant subrings of the enveloping algebras}

\bsno

Let $L$ be a finite dimensional Lie algebra over $k$ and $U(L)$ be
the universal enveloping algebra. There is a standard filtration given by
$F_i=(k+L)^i$ with the associated graded ring isomorphic to the
commutative polynomial ring $k[L]$. 

\proclaim{Lemma 6.1} If $g$ is a filtered automorphism of $U(L)$, there
there is a Lie algebra automorphism $\sigma: L\to L$ and a linear
map $\epsilon: L\to k$ with $\epsilon([L,L])=0$ such that 
$g|_{L}=\sigma+\epsilon$. Also $\hdet_{U(L)} g=\det \sigma|_L$.
\endproclaim

\demo{Proof} Since $g$ sends $L$ to $L\oplus k$, $g$ decomposes as 
$\sigma+\epsilon$. For every $a, b\in L$, 
$$\sigma([a,b])=g([a,b])-\epsilon([a,b])
=[\sigma(a)+\epsilon(a),\sigma(b)+\epsilon(b)]-\epsilon([a,b])
=[\sigma(a),\sigma(b)]-\epsilon([a,b]).$$
Therefore $\sigma$ is an endomorphism (and hence automorphism) of 
Lie algebra $L$ and $\epsilon([L,L])=0$. By the definition and 
Lemma 2.1(1), $\hdet_{U(L)} g=\det \sigma|_L$. \hfill $\square$
\enddemo

In the rest of this section we study the cases when the determinant
is 1. We assume that $\charr k=0$. The next lemma is elementary.

\proclaim{Lemma 6.2} Let $\nu$ be a nilpotent linear map of 
a finite dimensional space $V$, then the linear map 
$$\exp(\nu):=\sum_{n=0}^{\infty}\frac{\nu^n}{n!}$$ is a linear map of $V$
with determinant $1$.
\endproclaim

The above lemma is also true if the trace of $\nu$ is zero.
An automorphism $\sigma$ of the Lie algebra $L$ is called {\it inner} if
$\sigma=\exp(\ad x)$ for some $x\in L$ and $\ad x$ is nilpotent. 
The {\it inner automorphism group} of $L$, denoted by $\Int L$, is
generated by all inner automorphisms of $L$. By [{\bf Ja}, Chapter IX],
$\Int L$ is a normal subgroup of the automorphism group $\Aut(L)$ of $L$.

\proclaim{Corollary 6.3} Let $L$ be a Lie algebra over $k$ with 
$\charr k=0$.

(1) If $g|_L\in \Int L$, then $\hdet_{U(L)} g=1$.

(2) If $G$ is a finite subgroup of $\Int L$, then $U(L)^G$ is
Auslander-Gorenstein and $\GKdim$-Macaulay.
\endproclaim 

\demo{Proof} (1) A consequence of Lemmas 6.1 and 6.2.

(2) A consequence of (1) and Theorem 3.5(3).
\hfill $\square$
\enddemo

The following classification of automorphism groups of simple Lie algebra
is standard (see [{\bf Ja}], [{\bf Hu}], [{\bf Ka}]).
Let $k$ is an algebraically closed field of characteristic $0$. Let
$\Gamma L=\Aut L/\Int L$ [{\bf Hu}, p. 66].  Then

(a) $\Gamma L$ is the graph automorphism of the Dynkin diagram. 

(b) $\Gamma L=1$ (namely, $\Aut L=\Int L$) for the following types: 
$B_n, C_n$, $E_7, E_8$, $F_4$, $G_2$.

(c) $\Gamma$ of other types are

$\qquad$ $\Gamma A_n={\Bbb Z}/2{\Bbb Z}$ $(n\geq 2)$
and the outer automorphism  
is given by
$$g:\quad e_i\to e_{n+1-i}, \quad f_i\to f_{n+1-i}, \qquad i=1, \cdots, n $$

$\qquad$ $\Gamma D_4=S_3=<x, y>$ where 
$$\align x:&\quad e_1\to e_3,\quad e_3\to e_4,\quad e_4\to e_1,\quad e_2\to e_2\\
y:&\quad e_3\to e_4,\quad e_4\to e_3,\quad e_1\to e_1,\quad e_2\to e_2.
\endalign
$$

$\qquad$ $\Gamma D_l={\Bbb Z}/2{\Bbb Z}$ $(n> 4)$ 
and the outer automorphism is given by
 $$g:\quad e_{n-1}\to
e_n, \quad e_{n}\to e_{n-1}, \quad e_i\to e_i.$$ 

$\qquad$ $\Gamma E_6={\Bbb Z}/2{\Bbb Z}$ and the outer
automorphism is given by
$$g:\quad e_1\to e_6,\quad e_2\to e_5, \quad e_2\to e_2,\quad e_4\to e_4.$$

In the above we assume that $\{e_i, f_i, h_i\}$ are the Serre generators of the simple Lie algebra. We assume the similar formulas for the action of the outer automorphism
on $f_i$ and $h_i$. 


\proclaim{Corollary 6.4} 
Let $L$ be a simple Lie algebra
of the following types over an algebraically closed field $k$ with $\charr k=0$, and $G$ is a finite group of filtered automorphisms of $U(L)$:
 
(a) $A_n$, $n=4m$ or $n=4m+1$;

(b) $B_n$;

(c) $C_n$;

(d) $E_6$, $E_7$, $E_8$,  $F_4$ and $G_2$. 

\noindent
 then 
$U(L)^G$ is Auslander-Gorenstein and $\GKdim$-Macaulay. 
\endproclaim

\demo{Proof} By Corollary 6.3 and the above description of
$\Aut(L)$ it suffices to relevant outer automorphisms
are of determinant $1$.

Let $g$ be the outer automorphism given above (in case
of $D_4$, $g
=x $ or $y$). We compute directly that

(1) $\det g=(-1)^{n(n+3)/3}$   for $A_n$; 

(2) $\det x=1, \det y=-1$ for $D_4$;

(3) $\det g=-1$ for $D_n$;

(4) $\det g=1$ for $E_6$.

\noindent
Then the result follows. \hfill $\square$
\enddemo

\noindent
{\bf Remark 6.5.} In the case of $D_4$, if $G/(\Int(L)\cap G) \subset
C_3=\langle x\rangle$, then the result still holds.

\medskip

Theorem 0.5 follows from Corollaries 6.3 and 6.4.

\msno
{\bf Remark 6.6.} Alev and Polo showed that if $L$ is a semisimple
Lie algebra over an algebraically closed field of characteristic zero
and if $G$ is a non-trivial finite group of automorphisms of $U(L)$, then 
$U(L)^G$ is not isomorphic to any enveloping algebra $U(L')$ [{\bf AP}, Theorem 1]. 
Similar result holds for Weyl algebras [{\bf AP}, Theorem 2]. 

\msno
{\bf Remark 6.7.} It was informed to us by Lance Small that their
result [{\bf KS}, Proposition 4(5)] needs some extra hypothesis such as
the automorphisms of the Lie algebra are inner. Their method is similar to
 ours, namely, lifting the Gorenstein property from the associated graded
ring to the universal enveloping algebra. Note that the invariant subring 
of the associated graded ring needs not be Gorenstein if the 
automorphisms are not inner as the next example shows.

\msno{\bf Example 6.8.} Let $L=sl_3$ and $g$ is the outer automorphism
described before Corollary 6.4. Let $G=\{id, g\}$. Then 
$\G U(L)^G$ is not Gorenstein. In fact we compute that 
$$H_{\G U(L)^G}(t)=\frac12(\frac 1{(1-t^2)^3}+
\frac 1{(1-t)^6})=\frac{1+3t^2}{(1-t)^3(1-t^2)^3}.
$$
There is no integer $m$ such that $H_{\G U(L)^G}(t^{-1})
=\pm t^m H_{\G U(L)^G}(t)$, so $\G U(L)^G$ is not Gorenstein
([{\bf Ei}, Ex. 21.17(c)]). 



\bsno
\centerline{\bf Acknowledgment}

\bsno

The first author was supported by the NSA and the second author
was supported by the NSF and a Sloan Research Fellowship.

\bsno

\bsno

\centerline{\bf References}

\bsno

\item{[{\bf AHV}]} J. Alev, T. J. Hodges and J.-D. Velez, Fixed rings of the 
             Weyl algebra $A_1({\Bbb C})$, J. Algebra {\bf 130} (1990), 
             no. 1, 83--96.
\item{[{\bf AP}]} J. Alev and P. Polo, A rigidity theorem for finite group 
          actions on enveloping algebras of semisimple Lie algebras,
           Adv. Math. {\bf 111} (1995), no. 2, 208--226.
\item{[{\bf AST}]} M. Artin, W. Schelter and J. Tate, Quantum deformations of 
GL$_n$, Comm. Pure Appl. Math. {\bf 44} (1991), no. 8-9, 879--895.

\item{[{\bf ASZ}]} K. Ajitabh, S. P. Smith and J. J. Zhang,
        Auslander-Gorenstein rings, Comm. Algebra {\bf 26}
        (1998), no. 7, 2159--2180. 
\item{[{\bf Ben}]} D. J. Benson, Polynomial invariants of finite groups, 
         London Mathematical Society Lecture Note Series, 190, 
         Cambridge University Press, Cambridge, 1993.
\item{[{\bf B1}]} J.-E. Bj\"ork, The Auslander condition on Noetherian rings, 
           {\it  S\'eminaire d'Alg\`ebre Paul Dubreil et Marie-Paul 
            Malliavin}, 39\`eme Ann\'ee (Paris, 1987/1988), 137--173, 
            Lecture Notes in Math., 1404, Springer, Berlin-New York, 1989. 
\item{[{\bf B2}]} J.-E. Bj\"ork, Filtered Noetherian rings, {\it Noetherian rings 
            and their applications}, (Oberwolfach, 1983), 59--97, Math. 
            Surveys Monographs, 24, Amer. Math. Soc., Providence, RI, 1987.
\item{[{\bf BE}]} J.-E. Bj\"ork and E. K. Ekstr\"om, Filtered 
            Auslander-Gorenstein rings, {\it Operator algebras,
            unitary representations, enveloping algebras, and 
            invariant theory} (Paris, 1989), 425--448, Progr.
            Math., 92, Birkhäuser Boston, Boston, MA, 1990. 

\item{[{\bf Ei}]} D. Eisenbud, ``Commutative algebra, with a view toward 
            algebraic geometry'', Graduate Texts in Mathematics, 150,
            Springer-Verlag, New York, 1995.

\item{[{\bf Ek}]} E. K. Ekstr\"om, The Auslander condition on graded
        and filtered noetherian rings,
        {\it S\'eminaire d'Alg\`ebre Paul Dubreil et Marie-Paul Malliavin, 
        39\`eme Ann\'ee} (Paris, 1987/1988), 220--245, Lecture Notes 
        in Math., 1404, Springer, Berlin-New York, 1989. 
\item{[{\bf FKK}]} H. Fujita, E. Kirkman and J. Kuzmanovich, Global and Krull
           dimensions of quantum Weyl algebras, J. Algebra, to appear.
\item{[{\bf GZ}]} A. Giaquinto and J. J. Zhang, Quantum Weyl algebras, 
           J. Algebra {\bf 176} (1995), no. 3, 861--881. 

\item{[{\bf Hu}]} J. E. Humphreys, Introduction to Lie algebras and 
         representation theory, Second printing, revised,
         Graduate Texts in Mathematics, 9. Springer-Verlag, 
         New York-Berlin, 1978. 
\item{[{\bf Ja}]} N. Jacobson, ``Lie algebras'', Interscience Tracts in Pure 
          and Applied Mathematics, No. 10, Interscience Publishers 
          (a division of John Wiley \& Sons), New York-London, 1962 
\item{[{\bf JiZ}]} N. Jing and J. J. Zhang, On the trace of graded automorphisms,
         J. Algebra {\bf 189} (1997), no. 2, 353--376. 
\item{[{\bf JoZ}]} P. J{\o}rgensen and J. J. Zhang, Gourment's Guide to
             Gorensteinness, Adv. in Math., to appear.
\item{[{\bf Ka}]} V. G. Kac, Infinite dimensional Lie algebras, 3rd ed., Cambridge
Univ. Press, 1990.
\item{[{\bf KS}]} H. Kraft and L. W. Small, Invariant algebras and completely
          reducible representations, Math. Research Letters {\bf 1}
          (1994), 297-307. 
\item{[{\bf L1}]} T. Levasseur, Grade des modules sur certains anneaux filtr\'es 
         (French), Comm. Algebra {\bf 9} (1981), no. 15, 1519--1532. 
\item{[{\bf L2}]} T. Levasseur, Some properties of noncommutative
        regular rings, Glasgow Math. J. {\bf 34} (1992), 277-300.
\item{[{\bf Mo}]} S. Montgomery, {\it Fixed rings of finite automorphism 
           groups of associative rings,} Lecture Notes in Mathematics, 
           818. Springer, Berlin, 1980.
\item{[{\bf SS}]} S. P. Smith and J. M. Staniszkis, Irreducible 
         representations of the $4$-dimensional Sklyanin algebra at 
         points of infinite order. J. Algebra {\bf 160} (1993), 
         no. 1, 57--86. 
\item{[{\bf YZ}]} A. Yekutieli and J. J. Zhang, Rings with Auslander dualizing 
        complexes, J. Algebra, to appear.
\item{[{\bf Zh}]} J. J. Zhang, Connected graded Gorenstein algebras
        with enough normal elements, J. Algebra {\bf 189} (1997),
        390-405.

\enddocument